\documentstyle{amsppt}
\magnification=1200
\hsize=150truemm
\vsize=224.4truemm
\hoffset=4.8truemm
\voffset=12truemm

\NoRunningHeads

\define\C{{\bold C}}
 
\define\R{{\bold R}}
 
\let\thm\proclaim
\let\fthm\endproclaim

 \define\p{ P^1(\C) }
\define\pp{ P^2(\C) }
 \define\cs{  \C^*  }

\define\hyp{hyperbolicit\'e }
  \define\br{courbe de Brody}
 
\newcount\tagno
\newcount\secno
\newcount\subsecno
\newcount\stno
\global\subsecno=1
\global\tagno=0
\define\ntag{\global\advance\tagno by 1\tag{\the\tagno}}

\define\sta{\ 
{\the\secno}.\the\stno
\global\advance\stno by 1}

\define\stas{\the\stno
\global\advance\stno by 1}

\define\sect{\global\advance\secno by 1
\global\subsecno=1\global\stno=1\
{\the\secno}. }

\def\nom#1{\edef#1{{\the\secno}.\the\stno}}
\def\inom#1{\edef#1{\the\stno}}
\def\eqnom#1{\edef#1{(\the\tagno)}}

\newcount\refno
\global\refno=0

\def\nextref#1{
      \global\advance\refno by 1
      \xdef#1{\the\refno}}

\def\bref {\ref\global\advance\refno by 1\key{\the\refno}}


 \nextref\AUD
\nextref\BER 
\nextref\BRO
\nextref\DEB
\nextref\GRE
\nextref\GRO
\nextref\KOB
\nextref\KRU
\nextref\LEH
\nextref\SIK

\topmatter

\title 
Un th\'eor\`eme de Green presque complexe
 \endtitle

\author  Julien Duval
\endauthor

\abstract\nofrills{\smc R\'esum\'e. }\ \ On montre l'\hyp  du 
compl\'ementaire de 
cinq droites en position g\'en\'erale dans un plan projectif presque 
complexe, r\'epondant ainsi \`a une question de S. Ivashkovich. 

\null
{\smc \hskip 1,5cm  An almost complex version of a theorem by Green }

\null
\noindent
{\smc Abstract.}\ We prove the hyperbolicity 
of the complement of five lines in general position in an almost 
complex projective plane, answering a question by S. Ivashkovich.
 \endabstract

\endtopmatter 

\document

\subhead 0. Introduction \endsubhead
\stno=1

\null

Soit $X$ une vari\'et\'e de dimension $2n$ munie d'une structure 
presque complexe, i.e. d'un automorphisme $J$ de $TX$ tel que 
$J^2=-Id$. Quand $J$ n'est pas int\'egrable (non
localement \'equivalente \`a la structure complexe $i$ de $\C^n$), la 
vari\'et\'e
$X$ manque en g\'en\'eral d'objets holomorphes. Elle poss\`ede 
cependant des courbes J-holomorphes, i.e. des surfaces dont le plan tangent 
en tout point est 
une droite complexe pour $J$. 

En particulier elle a beaucoup de J-disques non constants
(voir l'article de J.-C. Sikorav dans [\AUD]). Un {\it J-disque} est 
une application $f:(D,i)\rightarrow (X,J)$ du disque unit\'e $D$ 
de $\C$ vers $X$, 
qui est
 J-holomorphe : elle v\'erifie $df\circ i=J\circ df$.
A la suite de Kruglikov et Overholt [\KRU], on d\'efinit donc
la pseudom\'etrique 
de Kobayashi-Royden $K$ sur $TX$ comme en complexe  
par :  
 $$K(p,v)= \text{Inf}\{1/r>0 /\text{il existe un J-disque avec }
f(0)=p 
\text{ et }d_0f( \partial /\partial x)=rv \},$$ 
o\`u $p$ est un point de $X$ et $v$ un 
vecteur de $T_p(X)$.

  La vari\'et\'e $X$ est 
dite {\it  hyperbolique} lorsque $K$ est une vraie m\'etrique. Au 
contraire, sa 
d\'eg\'en\'erescence se traduit par l'existence de J-disques 
 arbitrairement grands dans $X$ passant par un point dans une 
direction donn\'ee. On s'attend alors, au moins 
quand $X$ est compacte, \`a la 
pr\'esence de courbes enti\`eres dans la vari\'et\'e.

Pr\'ecis\'ement, appelons 
{\it \br} \ une application $f: (\C,i) \rightarrow (X,J)$ non 
constante,
 J-holomorphe et de d\'eriv\'ee born\'ee.
  
 Comme en complexe (voir Brody [\BRO]),   
 l'\hyp se caract\'erise ainsi [\KRU] :
  
\thm {Crit\`ere} Soit $(X,J)$ une vari\'et\'e presque complexe 
compacte. Alors $X$ est hyperbolique si et seulement si elle ne 
contient pas de \br.
\fthm

 Dans le cas complexe, on \'etudie plus  ais\'ement l'\hyp des
 compl\'ementai- res de diviseurs que celle des vari\'et\'es compactes. Le 
 crit\`ere pr\'ec\'edent demeure souvent valide. Ainsi, l'exemple de 
 base d\^u \`a Green [\GRE], l'\hyp du compl\'emen- taire de $2n+1$ 
 hyperplans en position g\'en\'erale dans $P^n(\C)$, se r\'eduit \`a 
 l'absence de \br \ dans $P^n(\C)$ \'evitant ces hyperplans. Il remonte 
 donc essentiellement \`a Picard en dimension complexe 1 et \`a 
 Borel en dimension sup\'erieure.

 A la suite de S. Ivashkovich, il est tentant 
 d'explorer cet exemple en presque complexe. Ceci n'a de sens 
 qu'en dimension r\'eelle 4 : en effet toute structure presque complexe est 
 int\'egrable en dimension 2, alors que l'on manque d'hypersur- faces
  J-holomorphes en dimension sup\'erieure \`a 4.  
 
 \null

 On se donne donc un {\it plan projectif presque complexe}. Autrement 
 dit, fixons une structure
 presque complexe $J$ 
 sur $\pp$  positive par rapport \`a la forme de 
 Fubini-Study $\omega$ : $\omega(\cdot ,J\cdot )>0$.   
 Appelons {\it J-droite} de notre plan l'analogue presque complexe d'une 
 droite projective, donc une courbe J-holomorphe plong\'ee dans 
 $(\pp ,J)$, diff\'eomorphe \`a $\p$ et de degr\'e 1 en 
 homologie.  
  
 D'apr\`es Gromov [\GRO] (voir aussi [\SIK]), un tel plan presque 
 complexe poss\`ede beaucoup de J-droites. Ainsi l'espace de ces J-droites 
 est diff\'eomorphe \`a $\pp$. Elles v\'erifient de plus les relations 
 d'incidence usuelles : par deux points 
 distincts passe une unique J-droite; deux 
 J-droites distinctes se coupent transversalement en un point unique;
 les J-droites passant par un point $p$ forment un pinceau diff\'eomorphe \`a 
 $\p$, donnant une projection centrale $\pi : \pp \setminus \{p\} 
 \rightarrow \p$.  
  
 Soit $C$ une r\'eunion de cinq J-droites en position 
 g\'en\'erale (i.e. sans point triple) de ce plan presque complexe. 
 Voici notre r\'esultat :
 
 \thm {Th\'eor\`eme} Le plan projectif presque complexe $(\pp ,J)$ ne 
 contient pas de \br  \ \'evitant la configuration $C$.
 \fthm
 \thm {Corollaire} Le compl\'ementaire $(\pp \setminus C,J)$ est 
 hyperbolique.
 \fthm
  Dans cette direction, Debalme et Ivashkovich [\DEB] avaient 
  auparavant
  remarqu\'e que
 l'\hyp de $(\pp \setminus C,J)$ \'etait une propri\'et\'e ouverte dans 
 l'espace des 
 configurations $(C,J)$.
  
 \null
 
 Le th\'eor\`eme se montre par l'absurde. Consid\'erons une courbe de 
 Brody $f$ 
  \'evitant la configuration $C$.  La remarque principale est que 
  son adh\'erence $F$ dans $\pp$ est
  "feuillet\'ee" 
 par des limites de $f$. 
 Par 
 positivit\'e d'intersection,  une telle feuille  
 \'evite
 une droite de $C$ ou y est contenue.  Or, par un analogue du 
 th\'eor\`eme de Liouville, $F$ doit rencontrer chaque droite de
 la configuration. 
Donc $F$ contient $C$ par propagation le long des feuilles.
 La contradiction est atteinte \`a un point double de $C$ puisque la 
 feuille y passant ne peut \^etre contenue \`a la fois dans les deux 
 droites 
 correspondantes.
 Le feuilletage de $F$ s'obtient par un argument de famille normale 
 en analysant la courbe de Brody $f$ sous les 
 projections centrales $\pi$ \`a partir des points doubles de $C$. Le 
 point crucial est que $\pi \circ f : \C
  \rightarrow \cs $ est quasiconforme et d'ordre fini, donc 
  essentiellement
   un 
  rev\^etement.  
  
  \null
  
  Notons que ce sch\'ema g\'eom\'etrique est int\'eressant m\^eme 
  dans le cas 
  complexe. Il r\'eduit le th\'eor\`eme de Green \`a un  
  fait analytique \'el\'ementaire de th\'eorie de distribution des valeurs : 
  une fonction 
  enti\`ere ne s'annulant pas et 
  d'ordre fini est une exponentielle de polyn\^ome (comparer avec [\BER]).  
  
  \null
  
  Afin de pr\'eciser ceci, d\'ebutons
  par des pr\'eliminaires sur les suites de J-disques et les 
  projections centrales dans un plan projectif presque complexe 
  $(\pp,J)$.

  \subhead 1. Pr\'eliminaires \endsubhead

  \null
   
  Les objets consid\'er\'es dans la suite sont de classe $C^\infty$ sauf 
  mention du 
  contraire, les convergences de suites d'applications \'etant
  localement uniformes.  
  
  \null
  
  a) {\bf Positivit\'e d'intersection.}  
  
    \null
    
  Soient deux J-disques non constants $f,g: (D,i) \rightarrow (\pp ,J)$. 
  Supposons-les d'images distinctes et s'intersectant en $f(0)=g(0)=p$.
  Alors cette intersection est isol\'ee et
  l'intersection homologique
  des deux disques en $p$  est strictement positive ([\GRO], voir 
  aussi l'article de D. McDuff 
  dans [\AUD]).  Ceci entra\^\i ne le :
  \thm {Fait} Soit $(f_n)$ une suite de J-disques 
  convergeant vers un J-disque $f$. On suppose que $f_n(D)$ \'evite une 
  J-droite $L$ du 
  plan presque complexe. Alors $f(D)$ \'evite encore $L$ ou y est 
  contenu.
   
  \fthm
  
   Sinon le disque $f(D)$ couperait positivement la 
   droite $L$. On trouverait donc une courbe ferm\'ee dans $f(D)$
   enla\c cant localement $L$. Ce serait encore le cas pour $f_n(D)$ 
   pour $n$ assez grand, contredisant le fait que $f_n(D)$ \'evite $L$.
   
   \null
   
   b) {\bf Familles normales.}
   
   \null
   
    Remarquons d'abord que, si une suite $(f_n)$ de J-disques 
    converge vers une application continue $f: D \rightarrow \pp$, 
   alors celle-ci est de classe $C^\infty$ et on a
   convergence des d\'eriv\'ees 
   successives. La limite $f$ est donc un J-disque. Ceci r\'esulte de 
   la r\'egularit\'e 
   elliptique de l'\'equation des courbes J-holomorphes (voir 
   l'article de J.-C. Sikorav dans 
   [\AUD]).
   
   Une suite de J-disques est {\it normale} si de toute sous-suite on peut
   extraire une suite convergente. La remarque pr\'ec\'edente et 
   le th\'eor\`eme d'Ascoli donnent le :
    
   \thm {Crit\`ere} Une suite $(f_n)$ de J-disques est normale si et 
   seulement si la suite $(\Vert df_n  \Vert)$ des normes de ses d\'eriv\'ees 
   est localement uniform\'ement born\'ee.  
  \fthm
  
  Ici  $\Vert df_n  \Vert$ est mesur\'ee dans les m\'etriques standard 
  de $\C$ et $\pp$.
  
  Une suite non normale de J-disques produit, quant \`a elle, une 
  \br \  par reparam\'etrage (cf. [\BRO] 
  et [\KRU]) :
  
  \thm {Lemme de Brody} Soit $(f_n)$ une suite non normale de J-disques. 
  Alors il existe une suite de contractions affines $(r_n)$ 
  de $\C$ convergeant vers un point de $D$ telle que $(f_n\circ 
  r_n)$ converge vers une \br \ apr\`es extraction.

  \fthm
 
  Rappelons que  cette \br \  est une application non constante $f:(\C 
  ,i) \rightarrow
 (\pp ,J)$, J-holomorphe et de d\'eriv\'ee 
  born\'ee  
  ($\Vert 
  df\Vert  \leq $constante).
  
  \null
  
  A ce stade, voyons comment le corollaire d\'ecoule du 
  th\'eor\`eme.

  Soit $C$ une configuration de cinq J-droites en 
  position g\'en\'erale dans $\pp$. Si $(\pp \setminus C,J)$ n'est pas 
  hyperbolique, 
  on obtient une suite de J-disques dont les d\'eriv\'ees  en l'origine 
  explosent, donc une suite non normale \'evitant $C$.
  Celle-ci
  produit une \br \  dans $(\pp,J)$. Elle doit
  \'eviter $C$, contredisant le th\'eor\`eme : sinon  
  cette courbe de Brody rencontrerait une droite de la configuration $C$;
   par le a) elle y serait contenue tout en \'evitant les quatre autres 
   droites;  
   on obtiendrait ainsi une courbe enti\`ere non constante dans 
   $\p \setminus$\{4\ points\}, ce qui est impossible par le 
   th\'eor\`eme de Picard. $\square$
  
   \null
   
  c) {\bf Eclatement presque complexe.}
  
  \null
  
  Soit $p$ un point de notre plan presque complexe. D'apr\`es [\GRO] 
  (voir aussi [\SIK]) passe par $p$ un pinceau de J-droites param\'etr\'e
  par $\p$, donnant une projection centrale $\pi : \pp \setminus \{p\} 
 \rightarrow \p$.  
  
 Redressons localement ce pinceau sur le pinceau lin\'eaire des 
 droites complexes de $\C^2$ en $0$. On construit pour cela
 un diff\'eomorphisme
 $\Phi$ pr\`es de $p$, en projetant chaque J-droite $L$ du pinceau 
 sur sa tangente $T_pL$ parall\`element \`a $T_pL^{\bot}$. 
 Ce diff\'eomorphisme est de classe $C^\infty$ hors de $p$ mais
 seulement $C^{1+Lip}$ en $p$.
 La structure presque complexe transport\'ee par $\Phi$, encore 
 not\'ee $J$, sera donc de classe $C^\infty$ hors de $0$ et
 Lipschitz en $0$. On 
 peut toujours supposer que $J(0) = i$. 
 
 D\'efinissons alors l'{\it \'eclat\'e} $X$ de $(\pp,J)$ en $p$ comme 
 l'\'eclat\'e complexe usuel $\widetilde {\C^2}$ en $0$ via $\Phi$.
 
 Par construction, la projection centrale $\pi$ se rel\`eve en une 
 fibration $\tilde \pi: X \rightarrow \p$. La structure $J$ donne par 
 rel\`evement une structure presque complexe $\tilde J$ sur $X 
 \setminus E$ o\`u $E$ est le diviseur exceptionnel de $X$. Les fibres 
 de $\tilde \pi$ sont $ \tilde J$-holomorphes hors de $E$.

 \thm {Lemme} La structure presque complexe $\tilde J$ admet un 
 prolongement Lipschitz au diviseur exceptionnel $E$.
 
 \fthm
 \demo {D\'emonstration} On le v\'erifie via $\Phi$. Notons 
 $q$ la projection de $\widetilde {\C^2}$ sur $\C^2$. On a, dans une des 
 deux cartes 
 de l'\'eclat\'e usuel, $q(x,t) = (x,tx)$.
 Hors du diviseur exceptionnel $E=(x=0)$, la structure relev\'ee 
 s'obtient par  $\tilde J = (dq)^{-1} \circ J \circ dq$. Les 
 horizontales $(t=$constante) \'etant $\tilde J$-holomorphes, la 
 structure $\tilde J$ est de la forme $$\tilde J = \pmatrix l & m \\ 0 & j 
 \endpmatrix$$ o\`u $j$, $l$ et $m$ sont, en chaque point, des $\R$-endomorphismes de 
 $\C$.
 En posant de la m\^eme mani\`ere $J=\pmatrix a & b \\ c & d 
 \endpmatrix$ et en explicitant la diff\'erentielle $dq$, 
 on obtient les formules suivantes pour $\tilde J$ au point $(x,t)$ :
 $$l =a\circ q +( b\circ q )t, \ \ 
  m =(b\circ q) x, \ \ 
  j =-x^{-1}t(b\circ q) x+x^{-1}(d\circ q) x.$$
 Comme $J(z)=i+O(\vert z \vert )$, on v\'erifie bien
 que $\tilde J (x,t) = i + O(\vert x\vert)$. $\square$
 \enddemo

 \null

 d) {\bf Projections quasiconformes.}

 \null

 Comme dans le paragraphe pr\'ec\'edent, on se fixe $\pi$ une projection 
 centrale associ\'ee au pinceau de J-droites en $p$. Celle-ci n'est 
 pas en g\'en\'eral holomorphe. Cependant elle reste 
 quasiconforme en restriction aux J-disques du plan presque complexe.
 
 Pr\'ecisons ceci. On suppose que $\pi$ envoie l'orientation transverse
 du pinceau venant de $J$ sur celle de $\p$.

 Soit $0\leq k<1$.  
 Une application $g : 
 D  \rightarrow \p$ est dite {\it k-quasiconforme} si $\Vert 
 \bar \partial g \Vert \leq k\Vert \partial g \Vert $. Ici $\partial 
 g$ et  $\bar \partial g $ d\'esignent respectivement les 
 composantes $\C$-lin\'eaire et $\C$-antilin\'eaire de la 
 d\'eriv\'ee de $g$ pour les structures complexes  
 de $D$ et $\p$.  
 Notons que $g$ reste quasiconforme (pour une autre constante) si on 
 la compose par un diff\'eomorphisme pr\'eservant l'orientation de 
 $\p$.
 
 \thm {Proposition} Il existe une constante $k <1$ telle 
 que $\pi \circ f$ soit k-quasiconforme pour tout J-disque $f$ \'evitant $p$. 
 \fthm
 \demo {D\'emonstration} Relevons le J-disque $f$ \`a $X$ l'\'eclat\'e 
 presque complexe en $p$, en un $\tilde J$-disque $\tilde f$. On veut 
 voir que $\tilde \pi \circ \tilde f$ est quasiconforme.
 L'\'enonc\'e \'etant de nature locale, on peut supposer le disque 
 $\tilde f(D)$ contenu dans un ouvert de carte $U$ de $X$ 
 trivialisant la fibration. L'ouvert $U$ est donc diff\'eomorphe au 
 bidisque $D\times D$,  $\tilde \pi$ devenant la 
 projection sur le deuxi\`eme facteur. La structure $\tilde J$
 dans cette carte est alors
 de la 
 forme $$\tilde J=\pmatrix l & m \\ 0 & j 
 \endpmatrix.$$ 
 Comme on peut toujours supposer que $\tilde J (0) =i$,    
 la structure $j$ sera une petite
 perturbation $i+\epsilon$ de la 
 structure complexe de $D$ quitte \`a r\'etr\'ecir $U$. Ceci n'utilise 
 que la continuit\'e de $\tilde J$ (voir c)). 
  
 Ainsi $\tilde \pi \circ \tilde f$ se lit dans la carte comme 
 la deuxi\`eme projection $g: D \rightarrow D$ du $\tilde J$-disque. 
 D'apr\`es la forme 
 de $\tilde 
 J$, elle satisfait l'\'equation $dg \circ i= j(\tilde f)\circ dg = i\circ dg + 
 \epsilon(\tilde f)\circ dg$. On en d\'eduit la quasiconformalit\'e de 
 $g$ puisque $\epsilon$ est petit. L'application $\tilde \pi \circ \tilde f$ est 
 donc $k_U$-quasiconforme pour une constante ne d\'ependant que de la carte. 
 En recouvrant $X$ par un nombre fini de tels 
 ouverts $U$, on peut prendre pour $k$ le maximum des constantes $k_U$ 
 en question. $\square$
 
 \enddemo

 \noindent
 { \bf Remarque}. D'apr\`es le dernier chapitre de la 
 monographie de Lehto et Virtanen [\LEH], une application 
 k-quasiconforme peut toujours s'\'ecrire 
 comme compos\'ee 
 $h\circ \phi$ d'une fonction holomorphe $h$ et d'un hom\'eomorphisme 
 quasiconforme $\phi$ gr\^ace au th\'eor\`eme d'Ahlfors-Bers.
 Nous renvoyons \`a [\LEH] pour les d\'efinitions analytique et 
 g\'eom\'etrique des {\it 
 hom\'eomorphismes}
 quasiconformes, et leurs propri\'et\'es. 
  
  \null
  
 On en d\'eduit d\'ej\`a cet analogue du th\'eor\`eme
 de Liouville :
 \thm {Corollaire} Soit $f: (\C,i) \rightarrow (\pp,J)$ une 
 courbe enti\`ere J-holomorphe. On suppose que l'adh\'erence de 
 son image \'evite une J-droite $L$. Alors $f$ est constante.
 \fthm
 \demo{D\'emonstration} Plongeons $L$ dans un pinceau de J-droites. On 
 suppose que la projection centrale $\pi$ envoie $L$ sur le point \`a 
 l'infini de $\p$. Notons $g$ la compos\'ee $\pi \circ f$.
 Comme $\overline {f(\C)}$ \'evite $L$,   
  $\overline {g(\C)}$ \'evite l'infini. Autrement dit $g:\C 
  \rightarrow \C$ est born\'ee. 
 Par la proposition et la remarque, 
 $g$ s'\'ecrit comme une compos\'ee $h\circ \phi$ o\`u 
 $h$ est une fonction enti\`ere born\'ee et $\phi$ un 
 hom\'eomorphisme de $\C$. Par le th\'eor\`eme de 
 Liouville, $h$ et donc $g$ sont constantes. 
 Ainsi la courbe $f(\C)$ est contenue dans une J-droite $L'$ du pinceau. 
Or l'adh\'erence 
 $\overline {f(\C)}$ \'evite le point de rencontre de $L$ et $L'$ 
 que l'on identifie au point \`a l'infini de $L'$. Donc $f$ est
 constante par une nouvelle application du th\'eor\`eme de Liouville.
 $\square$
 \enddemo

 Abordons maintenant la d\'emonstration proprement dite du th\'eor\`eme.
 
 \null
 
 \subhead 2. D\'emonstration \endsubhead
 
 \null
 
 Rappelons que l'on raisonne par l'absurde. Soit $f:(\C,i) 
 \rightarrow (\pp,J)$ une courbe de Brody \'evitant une 
 configuration $C=\cup  L_i$ de cinq J-droites en position
 g\'en\'erale. Notons 
 $F= \overline  {f(\C)}$  l'adh\'erence de son image dans $\pp$. On a 
 le :
 \thm {Lemme g\'eom\'etrique} Le compact $F$ est r\'eunion de 
 J-disques $\Delta$ non constants obtenus comme limites de disques de 
 $f(\C)$.
 \fthm
 
 Admettons provisoirement ce lemme. Voici comment en d\'ecoule la 
 contradiction. Remarquons d\'ej\`a que les disques 
 $\Delta$ satisfont l'alternative suivante vis-\`a-vis de $C$ (cf. 
  1.a)) : ou $\Delta$ \'evite $C$, ou $\Delta$ est contenu dans l'une 
 des droites de $C$ en \'evitant les autres. En particulier, $F$ 
 doit \'eviter les points doubles de la configuration $C$. Cependant, 
 par l'analogue du th\'eor\`eme de Liouville (cf. 1.d)), $F$ rencontre 
 chaque droite $L_i$ de $C$. De plus l'intersection $F \cap L_i$, qui est ferm\'ee
 dans 
 $L_i$, y est aussi ouverte : le disque $\Delta$ passant par un point de 
 $F \cap L_i$ 
 est n\'ecessairement contenu dans $L_i$. Donc $F$ contient 
 toute la configuration. Contradiction. $\square$
 
 \demo {D\'emonstration du lemme g\'eom\'etrique} Elle repose sur 
 l'analyse de $F$ sous les projections centrales issues des points 
 doubles de la configuration. Soit $p$ un tel point double, par 
 exemple le point de rencontre de $L_1$ et $L_2$. Notons $\pi:\pp \setminus 
 \{p\} \rightarrow \p$ la projection correspondante. On suppose 
 que les droites $L_1$ et $L_2$ s'envoient par $\pi$ sur l'origine et 
 l'infini de $\p$. La 
 restriction de $\pi$ \`a la courbe de Brody $f$ est donc \`a valeurs 
 dans $\cs$. Le point crucial est que $\pi \circ f$ est presque un 
 rev\^etement :
  
 \thm{Lemme analytique} La compos\'ee  $g=\pi \circ f:\C \rightarrow 
 \cs$ est de la forme $e^{P\circ \phi}$ o\`u $P$ est un polyn\^ome 
 non constant et $\phi$ un hom\'eomorphisme de $\C$.
 \fthm
  En particulier, soit $\delta$ un disque
  dans $\cs$; alors  $g^{-1}(\delta)=
  U\cup \bigcup D_n$  o\`u $U$ et les disques $D_n$ sont disjoints,   
  $g\vert _U$ de degr\'e fini et $g\vert _{D_n}$ un 
  hom\'eomorphisme sur $\delta$.
   
  \null
  
  Admettons pour l'instant ce r\'esultat. Voici comment il 
  entra\^\i ne le lemme g\'eom\'etrique. 
  
  Soit $z$ un point de 
  $F\setminus f(\C)$. Quitte \`a permuter les droites, supposons 
  que $L_1$ et $L_2$ \'evitent $z$. Notons $t=\pi(z)$ et $\delta$ un 
  disque centr\'e en $t$ dans $\cs$. Par hypoth\`ese, $z$ est limite 
  d'une suite de points distincts $(z_n)$ de $f(\C)$. Par la remarque 
  suivant le lemme analytique on trouve, pour $n$ 
  assez grand, un disque
  $\Delta_n$ contenu 
  dans la courbe de Brody passant par $z_n$ tel que $\pi : 
  \Delta_n \rightarrow \delta $ soit un hom\'eomorphisme. 
  
  Il suffit 
  de voir que $(\Delta_n)$ forme une famille normale pour conclure.
  Pour cela reparam\'etrons conform\'ement les disques $\Delta_n$ par
  $f_n:(D,i) 
  \rightarrow (\pp,J)$ avec $\pi \circ f_n(0)=t$. Les 
  projections $g_n=\pi \circ f_n :D \rightarrow \delta$ sont ainsi
  des diff\'eomorphismes 
  k-quasiconformes (cf. 1.d)) envoyant l'origine sur $t$. Quitte \`a 
  extraire, on peut 
  donc supposer que $(g_n)$ converge vers un 
  hom\'eomorphisme $\gamma$ de $D$ sur $\delta$ [\LEH]. 
  
  On en d\'eduit bien que la suite de J-disques $(f_n)$ est normale. 
  Sinon, par le lemme de Brody (cf. 1.b)), on trouverait une suite de 
  contractions affines $(r_n)$ de $\C$ convergeant vers un point $a$ 
  de $D$ telle que $(f_n \circ r_n)$ converge vers une courbe de Brody 
  $h$ quitte \`a extraire.
  Ainsi $(g_n \circ r_n)$ aurait pour limite $\gamma(a)$. Autrement 
  dit, la courbe $h$ serait contenue dans la fibre de $\pi$ 
  au-dessus de $\gamma(a)$, donc dans une J-droite $L$ du 
  pinceau en $p$. Par ailleurs cette courbe de Brody \'eviterait encore
   $C$ comme limite d'une suite de J-disques hors de $C$ (cf. 1.b)). 
  Or $L \cap C$ contient au moins trois points. On obtiendrait donc une 
  courbe enti\`ere non constante dans $\p \setminus \{3$ points\}, 
  contredisant
  le th\'eor\`eme de Picard. $\square$
   
 \enddemo   
 
 \demo {D\'emonstration du lemme analytique} Elle repose sur le fait 
 que la courbe de Brody $f$ est d'ordre fini, ainsi que sa projection
 $g=\pi\circ f$. Comme $g$ est quasiconforme, on l'\'ecrira comme 
 compos\'ee $h\circ \phi$ d'une fonction enti\`ere $h$ 
  d'ordre fini et d'un hom\'eomorphisme. Or on sait bien qu'une 
  fonction enti\`ere ne 
  s'annulant pas et d'ordre fini est une exponentielle de polyn\^ome 
  (cf. par exemple [\BER]).
  
  D\'etaillons ceci. Notons $D_r$ le disque centr\'e en l'origine 
 de rayon $r$ dans $\C$. Rappelons que $\omega$ d\'esigne la forme 
 de Fubini-Study de $\pp$, celle de $\p$ \'etant not\'ee $\omega'$.  
 D\'efinissons, comme Nevanlinna et Ahlfors (voir [\KOB]),
 les {\it caract\'eristiques} de $f$ et $g$ par :
 $$T(f)(r)= \int_0^r ( \int_{D_\rho}f^*\omega )\tfrac {d\rho } 
 {\rho},\  \  \ 
 T(g)(r)= \int_0^r (\int_{D_\rho}g^*\omega')
 \tfrac {d\rho } {\rho}=
  \int_0^r (\int_{D_\rho}f^*\pi^*\omega')\tfrac {d\rho } {\rho}.$$
 Les applications $f$ ou $g$ sont {\it d'ordre fini} si leur caract\'eristique
 cro\^\i t au plus polynomialement 
 en $r$. C'est le cas pour $f$ : en effet sa d\'eriv\'ee est born\'ee, donc
 $f^*\omega$ est une 2-forme born\'ee sur $\C$ et $T(f)$ 
 cro\^\i t quadratiquement.  
 
 Pour passer \`a $g$, on va  
 comparer la forme singuli\`ere
 (le courant) $\pi^*\omega'$ \`a $\omega$. 
   Leur diff\'erence n'est plus tout \`a fait, comme en complexe, le $dd^c$ 
   d'un potentiel 
   \`a singularit\'e 
   logarithmique en $p$ :
    
  \thm {Fait} Il existe une fonction n\'egative $u$ et une 1-forme 
  born\'ee $\alpha$ sur $\pp \setminus \{p\}$ telles que 
  $\pi^*\omega'=\omega + dd^c_J u + d\alpha$ en dehors de $p$, 
  avec la notation $d^c_J u=-du\circ J$.
  \fthm

  En effet, comme $\omega$ et $\pi^*\omega'$ sont cohomologues, 
  il s'agit d'un probl\`eme local pr\`es de $p$ : 
  \'ecrire $\pi^*\omega'$
 comme somme du $dd^c_J$ d'un potentiel 
 n\'egatif et d'une diff\'erentielle d'une 1-forme born\'ee.  Or, 
 apr\`es redressement du pinceau en $p$ par 
 $\Phi$ (cf. 1.c), on a 
 $\pi^*\omega'=dd^clog\vert z \vert= dd^c_Jlog\vert z \vert + d\beta$
 o\`u $\beta =dlog\vert z \vert \circ (J-i)$. Cette 1-forme est bien
born\'ee pr\`es de $0$ car $J(z)-i = O(\vert z \vert)$  ($\Phi$ 
est $C^{1+Lip}$ en $p$). 

\null

 Puisque $f$ est J-holomorphe, on peut maintenant comparer 
 $T(g)$ et $T(f)$ :
      
  $$T(g)(r)= T(f)(r)+\int_0^r (\int_{D_\rho}dd^c(u\circ f)\tfrac {d\rho } {\rho} 
  +\int_0^r (\int_{\partial D_\rho}f^*\alpha)\tfrac {d\rho } 
  {\rho}.$$
   
 Classiquement [\KOB] la premi\`ere int\'egrale vaut
 $\int_0^{2\pi}u\circ f(re^{i\theta})d\theta - 2\pi u\circ f(0)$. Elle reste 
 donc major\'ee puisque $u$ est n\'egative. 
  La seconde cro\^\i t lin\'eairement : en effet,
  la 1-forme $f^*\alpha$ 
  est born\'ee sur $\C$ car la d\'eriv\'ee de $f$  et $\alpha$ le sont.  
  Donc $g$ est, comme $f$, d'ordre fini.
  
  \null
  
   Ecrivons maintenant $g$ comme compos\'ee
    $g=h\circ \phi$ d'une fonction enti\`ere $h$   
   et d'un hom\'eomorphisme quasiconforme $\phi$ de $\p$ fixant 
   l'infini (cf. 1.d) et le dernier chapitre de [\LEH]).   
   Notons que l'hom\'eomorphisme $\phi$ est H\"older pour la m\'etrique de 
   Fubini-Study de $\p$ [\LEH].  Puisqu'il fixe l'infini, son 
   inverse cro\^\i t polynomialement pour la m\'etrique usuelle de $\C$ : 
   on aura
   $\phi^{-1}(D_r) \subset  
  D_{r^d}$ pour un certain entier $d$ et $r$ assez grand.  
  Ceci entra\^\i ne que $h$ est, comme $g$, d'ordre fini :
  $T(h)(r)$ vaut en effet par d\'efinition 
   $$\int_0^r  \text 
   {Aire}_{\omega'}(g(\phi^{-1}(D_{\rho}))) \tfrac {d\rho } {\rho} 
     \leq   
    \int_0^r  \text 
   {Aire}_{\omega'}(g(D_{\rho^d})) \tfrac {d\rho } {\rho} +O(1)=  
   d^{-1}T(g)(r^d)+O(1).$$
  
   Rappelons enfin bri\`evement pourquoi $h$ 
   est une exponentielle de polyn\^ome.
   Comme $h$ ne s'annule pas sur $\C$, elle a un 
   logarithme $k$. Montrons que c'est un polyn\^ome.  
   Notons pour cela que, dans $\C$, $\omega' =dd^clog(1+\vert z \vert^2)$. 
   Donc :  
   $$T(h)(r)=\int_0^r (\int_{D_\rho}dd^clog(1+\vert h \vert^2)
   \tfrac {d\rho } {\rho} 
  =\int_0^{2\pi}
   log(\vert h \vert+\vert h \vert^{-1})(r e^{i\theta})
   d\theta +O(1).$$
   La deuxi\`eme \'egalit\'e r\'esulte de $(*)$ et du fait que $log\vert h\vert$ 
   est harmonique.
   Ainsi les moyennes sur les cercles de centre $0$ et de rayon $r$ de
   $\big\vert 
  log\vert h\vert \big\vert  =\vert Re(k)\vert$ ont une croissance 
  polynomiale. Or le d\'eveloppement en s\'erie enti\`ere de 
   $k$ \`a l'origine donne :
   $$ k^{(n)}(0)=\pi^{-1}n!\ r^{-n}\int_0^{2\pi}
   Re(k)(r e^{i\theta}) e^{-in\theta}
   d\theta.$$
   En faisant tendre $r$ vers l'infini, on obtient bien que 
   $k^{(n)}(0)=0$ pour $n$ grand. $\square$

 \enddemo
 \Refs

\widestnumber\no{99}
\refno=0
\bref \by M. Audin and J. Lafontaine ed. \book Holomorphic curves in 
symplectic geometry \bookinfo Progress in Math. \vol 117 \publ 
Birkh\"auser \yr 1994  \publaddr Basel
\endref
 
\bref \by F. Berteloot et J. Duval \paper Sur l'hyperbolicit\'e de 
certains compl\'ementaires \jour Ens. Math. 
\vol47\yr2001\pages253--267
\endref
 
\bref \by R. Brody \paper Compact manifolds and hyperbolicity \jour 
Trans. Amer. Math. Soc. \vol235\yr1978\pages213--219
\endref

\bref \by R. Debalme and S. Ivashkovich \paper Complete hyperbolic 
neighborhoods in almost-complex surfaces \jour Int. J. of Math. 
\vol12\yr2001\pages211--221
\endref

 \bref \by M. Green \paper Some Picard theorems for holomorphic maps 
 to algebraic varieties \jour Amer. J.  Math. 
 \vol97\yr1975\pages43--75
 \endref
 
 \bref \by M. Gromov \paper Pseudo holomorphic curves in symplectic 
 manifolds \jour Invent. Math. \vol83\yr1985\pages307--347
 \endref
 
 \bref \by S. Kobayashi \book Hyperbolic complex spaces \bookinfo 
 Grund. der math. Wiss.\vol 318 \publ 
 Springer \yr 1998 \publaddr Berlin 
 \endref
 
  \bref \by B. Kruglikov and M. Overholt \paper Pseudoholomorphic 
  mappings and Kobayashi hyperbolicity \jour Diff. Geom. Appl. 
  \vol11\yr1999\pages265--277
  \endref
  
  \bref \by O. Lehto and K.I. Virtanen \book Quasiconformal mappings 
  in the plane \bookinfo Grund. der math. Wiss. \vol 126 \publ 
  Springer \yr 1973 \publaddr Berlin
  \endref
  
  \bref \by J.-C. Sikorav \paper Dual elliptic planes \jour preprint 
  2000
  arXiv math.SG/0008234
  \endref

\endRefs

\address 
\noindent  
Laboratoire \'Emile Picard, UMR CNRS 5580, 
  Universit\'e Paul Sabatier, 31062 Toulouse Cedex 4, France.
 \endaddress
\email 
  duval\@picard.ups-tlse.fr 
\endemail

\enddocument